%% file: 3nekom.tex
\documentclass[11pt,reqno]{amsart}
\usepackage{amssymb,amscd,amsbsy}
\usepackage{amssymb,amscd,amsbsy,mathrsfs}
\input{classstartscr}

\theoremstyle{remark}

\newtheorem*{rem*}{Remark}

\newcommand\Li{{\rm Lip}}
\newcommand\fM{\frak M}

\newcommand\dg{\frak D}

\begin{document}

\newcommand{\vse}{\vspace{.2in}}
\numberwithin{equation}{section}

\title{Functions of triples of noncommuting self-adjoint operators and their perturbations}
\author{V.V. Peller}
\thanks{the author is partially supported by NSF grant DMS 1300924}

\begin{abstract}
In this paper we study properties of functions of triples of not necessarily commuting self-adjoint operators. The main result of the paper shows that unlike in the case of functions of pairs of self-adjoint operators there is no Lipschitz type estimates in the trace norm for arbitrary functions in the Besov class $B_{\be,1}^1(\R^3)$. In other words, we prove that there is no constant $C>0$ such that the inequality
\begin{align*}
\|f(A_1,B_1,C_1)&-f(A_2,B_2,C_2)\|_{\bS_1}\\[.1cm]
&\le C\|f\|_{B_{\be,1}^1}
\max\big\{\|A_1-A_2\|_{\bS_1},\|B_1-B_2\|_{\bS_1},\|C_1-C_2\|_{\bS_1}\big\}
\end{align*}
holds for arbitrary function $f$ in $B_{\be,1}^1(\R^3)$ and for arbitrary finite rank self-adjoint operators $A_1,\,B_1,\,C_1,\,A_2,\,B_2$ and $C_2$.
\end{abstract}

\maketitle


\

%
%
%
%

\setcounter{section}{0}
\section{\bf Introduction}
\setcounter{equation}{0}
\label{In}

\

The spectral theorem for commuting self-adjoint operators implies that for commuting 
self-adjoint operators $A_1$ and $A_2$ and for a Lipschitz function $f$ on the real line
$\R$ the following Lipschitz type estimate holds
$$
\|f(A_1)-f(A_2)\|\le\|f'\|_{L^\be(\R)}\|A_1-A_2\|.
$$
The same inequality holds for the norms in Schatten--von Neumann classes $\bS_p$ with
$p\ge1$. However, for noncommuting self-adjoint operators, the situation is quite different.
A Lipschitz function $f$ on $\R$ does not have to be {\it operator Lipschitz}, i.e., the inequality
$$
|f(x_1)-f(x_2)|\le\const|x_1-x_2|,\quad x_1,~x_2\in\R,
$$
does not imply that
$$
\|f(A_1)-f(A_2)\|\le\const\|A_1-A_2\|
$$
for self-adjoint operators $A_1$ and $A_2$.
This was proved by Farforovskaya in \cite{F1}. She also proved in \cite{F2} that
there exist a Lipschitz function $f$ on $\R$ and self-adjoint operators $A_1$ and $A_2$ such that $A_1-A_2$ belongs to trace class $\bS_1$, but $f(A_1)-f(A_2)\not\in\bS_1$.

It turns out that a function $f$ on $\R$ is operator Lipschitz if and only if it takes trace class perturbations to trace class increments, i.e.,
$$
A=A^*,\quad B=B^*,\quad A-B\in\bS_1\quad\Longrightarrow\quad f(A)-f(A)\in\bS_1
$$
if we consider not necessarily bounded self-adjoint operators $A$ and $B$, see \cite{AP}.

It was shown later in \cite{Mc} and \cite{Ka} that the function $x\mapsto|x|$ is not operator Lipschitz. Necessary conditions for operator Lipschitzness were obtained in \cite{Pe2} and \cite{Pe3}. In particular, it was proved in \cite{Pe2} that operator Lipschitz functions on $\R$ must belong locally to the Besov class $B_{1,1}^1(\R)$. Note that this was deduced from the trace class criterion for Hankel operators, see \cite{Pe1} and \cite{Pe4}.

On the other hand, it was proved in \cite{Pe2} and \cite{Pe3} that functions in the Besov class $B_{\be,1}^1(\R)$ are necessarily operator Lipschitz.
This result was generalized in \cite{APPS} to functions of normal operators. It was shown in \cite{APPS} that if $f$ is a function of two variables that belongs to the Besov class  
$B_{\be,1}^1(\R^2)$, then $f$ is an {\it operator Lipschitz function on $\R^2$}, i.e.,
$$
\|f(N_1)-f(N_2)\|\le\const\|f\|_{B_{\be,1}^1}\|N_1-N_2\|
$$
for arbitrary normal operators $N_1$ and $N_2$. The same Lipschitz type inequality holds in the Schatten--von Neumann norm $\bS_p$ for $p\ge1$.
Note also that in \cite{NP} this result was generalized to the case of functions of $d$-tuples of commuting self-adjoint operators: if $f$ belongs to the Besov class $B_{\be,1}^1(\R^d)$ and $(A_1,\cdots,A_d)$ and $(B_1,\cdots,B_d)$ are $d$-tuples of commuting self-adjoint operators, then
$$
\|f(A_1,\cdots,A_d)-f(B_1,\cdots,B_d)\|\le\const\|f\|_{B_{\be,1}^1}
\max_{1\le j\le d}\|A_j-B_j\|
$$
and the same inequality holds for Schatten--von Nemann norms $\bS_p$ with $p\ge1$.

Let me also mention that in \cite{KPSS} it was shown that for an arbitrary Lipschitz function $f$ on $\R^d$ and for $p\in(1,\be)$ the following Lipschitz type inequality holds:
$$
\|f(A_1,\cdots,A_d)-f(B_1,\cdots,B_d)\|_{\bS_p}\le\const\|f\|_{\Li}
\max_{1\le j\le d}\|A_j-B_j\|_{\bS_p}
$$
for arbitrary $d$-tuples of commuting self-adjoint operators $(A_1,\cdots,A_d)$ and $(B_1,\cdots,B_d)$. Earlier in the case $d=1$ this was established in \cite{PS}.

We refer the reader to the survey article \cite{AP}, which is a comprehensive study of operator Lipschitz functions. 

The problem of the behavior of functions of pairs of {\it noncommuting} self-adjoint operators under perturbation was studied in \cite{ANP}. For a pair $(A,B)$ of not necessarily commuting self-adjoint operators the functions $f(A,B)$ are defined as double operator integrals:
$$
f(A,B)=\iint f(x,y)\,dE_A(x)\,dE_B(y)
$$
under the assumption that the double operator integral makes sense. Here  $E_A$ and $E_B$ stand for the spectral measures of $A$ and $B$.

In the case when $A$ and $B$ are finite rank self-adjoint operators (or, more general, if $A$ and $B$ have finite spectra), the operator $f(A,B)$ is defined for all functions $f$ on $\R^2$:
$$
f(A,B)=\sum_{j,k}f(\l_j,\mu_k)P_jQ_k,
$$
where 
$$
A=\sum_j\l_jP_j\quad\mbox{and}\quad B=\sum_k\mu_kQ_k
$$
are the spectral expansions of $A$ and $B$.

It turned out that the situation in the case of noncommuting operators is different. It was shown in \cite{ANP} that if $f$ belongs to the Besov class $B_{\be,1}^1(\R^2)$ and $1\le p\le2$, then the following Lipschitz type estimate holds:
$$
\|f(A_1,B_1)-f(A_2,B_2)\|_{\bS_p}\le\const\|f\|_{B_{\be,1}^1}
\max\big\{\|A_1-A_2\|_{\bS_p},\|B_1-B_2\|_{\bS_p}\big\}
$$
for arbitrary pairs $(A_1,B_1)$ and $(A_2,B_2)$ of not necessarily commuting self-adjoint operators.

On the other hand, it was shown in \cite{ANP} that there is no Lipschitz type estimate in the norm of $\bS_p$ for $p>2$ as well as in the operator norm. In other words if $p>2$, there is no constant $C$ such that
$$
\big\|f(A_1,B_1)-f(A_2,B_2)\big\|_{\bS_p}\le C\|f\|_{B_{\be,1}^1}
\max\big\{\|A_1-A_2\|_{\bS_p},\|B_1-B_2\|_{\bS_p}\big\}
$$
for arbitrary finite rank self-adjoint operators $A_1$,$B_1$, $A_2$ and $B_2$.
The same is true in the operator norm.

In this paper we deal with functions of triples of not necessarily commuting self-adjoint operators. For a triple $(A,B,C)$ of not necessarily commuting self-adjoint operators and a function $f$ on $\R^3$, the operator $f(A,B,C)$ is defined as the triple operator integral
$$
f(A,B,C)=\iiint f(x,y,z)\,dE_A(x)\,dE_B(y)\,dE_C(z)
$$
in the case when the triple operator integral is defined. Again, if $A$, $B$ and $C$ have finite spectra, the triple operator integral on the right is well defined for all functions $f$ and
$$
f(A,B,C)=\sum_{\l\in\s(A),\,\mu\in\s(B),\,\nu\in\s(C)}f(\l,\mu,\nu)
E_A(\{\l\})E_B(\{\mu\})E_C(\{\nu\}).
$$

The main objective of this paper is to show that {\it unlike in the case of functions of two noncommuting self-adjoint operators}, there is no Lipschitz type estimate in the trace norm
for functions in the Besov class $B_{\be,1}^1(\R^3)$. In other words, there is no constant $C>0$ such that
\begin{align*}
\big\|f(A_1,B_1,C_1)&-f(A_2,B_2,C_2)\big\|_{\bS_1}\\[.2cm]
&\le C\|f\|_{B_{\be,1}^1}
\max\{\|A_1-A_2\|_{\bS_1},\|B_1-B_2\|_{\bS_1},\|C_1-C_2\|_{\bS_1}\}
\end{align*}
for arbitrary functions $f$ in $B_{\be,1}^1(\R^3)$ and arbitrary finite rank self-adjoint operators $A_1$, $B_1$, $C_1$, $A_2$, $B_2$ and $C_2$.

\

\section{\bf Multiple operator integrals}
\setcounter{equation}{0}
\label{dois}

\

Double operator integrals appeared in the paper \cite{DK} by Daletskii and S.G. Krein. 
Later the beautiful theory of double operator integrals was created by Birman and Solomyak in \cite{BS1}, \cite{BS2} and \cite{BS3}. 

Let $(\X,E_1)$ and $(\Y,E_2)$ be spaces with spectral measures $E_1$ and $E_2$
on a Hilbert space $\h$, let $T$ be a bounded linear operator on $\h$ and let $\Phi$ be a bounded measurable function on $\X\times\Y$. {\it Double operator integrals} are expressions of the form
\bay
\label{dvoi}
\int\limits_\X\int\limits_\Y\Phi(x,y)\,d E_1(x)T\,dE_2(y).
\ey
Birman and Solomyak's starting point is the case when $T$ belongs to the Hilbert--Schmidt class $\bS_2$. In this case they defined double operator integrals of the form \rf{dvoi}
for arbitrary bounded measurable $\Phi$ and proved that
$$
\left\|\int\limits_\X\int\limits_\Y\Phi(x,y)\,d E_1(x)T\,dE_2(y)\right\|_{\bS_2}
\le\|\Phi\|_{L^\be}\|T\|_{|bS_2}
$$
(see \cite{BS1}).

To define double operator integrals for arbitrary bounded linear operators $T$ in the general case, restrictions on $\Phi$ must be imposed. Double operator integrals for arbitrary bounded operators $T$ can be defined for functions $\Phi$ that are {\it Schur multipliers} with respect to the spectral measures $E_1$ and $E_2$, see \cite{BS1}, \cite{Pe2}, \cite{Pi} and \cite{AP} for details.

However, in this paper we need double operator integrals only in the case when the spectral measures $E_1$ and $E_2$ are atomic and have finitely many atoms. 
We say that a spectral measure $E$ on a set $\X$ {\it is atomic and has finitely many atoms} if all subsets of $\X$ are measurable and there are points $a_1,\cdots,a_n$ in $\X$, called the {\it atoms}, such that $E(\X\setminus\bigcup_{j=1}^na_j)=\0$ and $E(\{a_j\})\ne\0$, $1\le j\le n$.

In the case when the spectral measures $E_1$ and $E_2$ are atomic with finitely many atoms, we
can define double operator integrals of the form \rf{dvoi} for arbitrary functions $\Phi$ by
\bay
\label{konato}
\int\limits_\X\int\limits_\Y\Phi(x,y)\,d E_1(x)T\,dE_2(y)=
\sum_{j,k}\Phi(a_j,b_k)E_1(\{a_j\})TE_2(\{b_k\}),
\ey
where the $a_j$ and the $b_k$ are the atoms of $E_1$ and $E_2$.

Under these assumptions, the norm of the linear transformer
$$
T\mapsto\iint\Phi(x,y)\,d E_1(x)T\,dE_2(y)
$$
(both in the operator norm and in the trace norm) is equal to the norm of the matrix
$\{\Phi(a_j,b_k)\}$ in the space of matrix Schur multipliers, i.e., the norm of the matrix transformer
$$
\{\g_{jk}\}\mapsto\{\Phi(a_j,b_k)\g_{jk}\}
$$
in the operator norm (or in the trace norm), see \cite{AP}.

Double operator integrals play an important role in perturbation theory. In particular, a special role is played by the following formula:
\bay
\label{razdrai}
f(A)-f(B)=\iint\limits_{\R\times\R}\frac{f(x)-f(y)}{x-y}\,dE_A(x)(A-B)\,dE_B(y),
\ey
which holds for arbitrary self-adjoint operators $A$ and $B$ with bounded
$A-B$ and for arbitrary operator Lipschitz functions $f$ on $\R$, see \cite{BS3}
and \cite{AP}.

In this paper we consider only operators with finite spectra, in which case formula 
\rf{razdrai} holds for arbitrary functions $f$ on $\R$; moreover, the divided difference
$(x,y)\mapsto(f(x)-f(y))(x-y)^{-1}$ can be extended to the diagonal $\{(x,x):~x\in\R\}$
arbitrarily, i.e., the values of the divided difference on the diagonal do not affect the right-hand side of \rf{razdrai}. This can be verified elementarily.

{\it Multiple operator integrals}
$$
\underbrace{\int\cdots\int}_m\Phi(x_1,\cdots,x_m)
\,dE_1(x_1)T_1\,dE_2(x_2)T_2\cdots\,dE_{m-1}(x_{m-1})T_{m-1}\,dE_m(x_m)
$$
were defined for functions $\Phi$ in the (integral) {\it projective tensor product}
of $L^\be$ spaces in \cite{Pe5}. Later multiple operator integrals were defined in \cite{JTT} for functions $\Phi$ in the Haagerup tensor products of $L^\be$ spaces. We refer the reader to the survey article \cite{Pe6} for detailed information about multiple operator integrals.

Again, in this paper we consider only atomic spectral measures with finitely many atoms, in which case multiple operator integrals can be defined for arbitrary functions $\Phi$ by analogy with double operator integrals, see formula \rf{konato}.

\

\section{\bf Besov classes $\bs{B_{\be,1}^1(\R^d)}$}
\setcounter{equation}{0}

\

In this paper we need only Besov classes $B_{\be,1}^1(\R^d)$ of functions on the Euclidean space $\R^d$. We give here a brief introduction to such spaces and we refer the reader to \cite{Pee} for detailed information about Besov classes.

Let $w$ be an infinitely differentiable function on $\R$ such
that
\bay
\label{w}
w\ge0,\quad\supp w\subset\left[\frac12,2\right],\quad\mbox{and} \quad w(s)=1-w\left(\frac s2\right)\quad\mbox{for}\quad s\in[1,2].
\ey

We define the functions $W_n$, $n\in\Z$, on $\R^d$ by 
$$
\big(\F W_n\big)(x)=w\left(\frac{\|x\|_2}{2^n}\right),\quad n\in\Z, \quad x=(x_1,\cdots,x_d),
\quad\|x\|_2\df\left(\sum_{j=1}^dx_j^2\right)^{1/2},
$$
where $\F$ is the {\it Fourier transform} defined on $L^1\big(\R^d\big)$ by
$$
\big(\F f\big)(t)=\!\int\limits_{\R^d} f(x)e^{-{\rm i}(x,t)}\,dx,\!\quad 
x=(x_1,\cdots,x_d),
\quad t=(t_1,\cdots,t_d), \!\quad(x,t)\df \sum_{j=1}^dx_jt_j.
$$
Clearly,
$$
\sum_{n\in\Z}(\F W_n)(t)=1,\quad t\in\R^d\setminus\{0\}.
$$

With each tempered distribution $f\in{\mathscr S}^\prime\big(\R^d\big)$, we
associate the sequence $\{f_n\}_{n\in\Z}$,
\bay
\label{fn}
f_n\df f*W_n.
\ey
The formal series
$
\sum_{n\in\Z}f_n
$
is a Littlewood--Paley type expansion of $f$. This series does not necessarily converge to $f$.

Initially we define the (homogeneous) Besov class $\dot B^1_{\be,1}\big(\R^d\big)$ as the space of all
\lb$f\in{\mathscr S}^\prime(\R^n)$
such that
\bay
\label{Wn}
\{2^{n}\|f_n\|_{L^\be}\}_{n\in\Z}\in\ell^1(\Z)
\ey
and put
$$
\|f\|_{B^1_{\be,1}}\df\big\|\{2^{n}\|f_n\|_{L^\be}\}_{n\in\Z}\big\|_{\ell^1(\Z)}.
$$
According to this definition, the space $\dot B^1_{\be,1}(\R^n)$ contains all polynomials
and all polynomials $f$ satisfy the equality $\|f\|_{B^s_{p,q}}=0$. Moreover, the distribution $f$ is determined by the sequence $\{f_n\}_{n\in\Z}$
uniquely up to a polynomial. It is easy to see that the series 
$\sum_{n\ge0}f_n$ converges in ${\mathscr S}^\prime(\R^d)$.
However, the series $\sum_{n<0}f_n$ can diverge in general. It can easily be proved that the series
\bay
\label{ryad}
\sum_{n<0}\frac{\partial f_n}{\partial x_j},\quad \mbox{where}\quad 1\le j\le d,
\ey
converges uniformly on $\R^d$.

Now we can define the modified (homogeneous) Besov class $B^1_{\be,1}\big(\R^d\big)$. We say that a distribution $f$
belongs to $B^1_{\be,1}(\R^d)$ if \rf{Wn} holds and
$$
\frac{\partial f}{\partial x_j}
=\sum_{n\in\Z}\frac{\partial f_n}{\partial x_j},\quad
1\le j\le d,
$$
in the space ${\mathscr S}^\prime\big(\R^d\big)$ (equipped with the weak-$*$ topology). Now the function $f$ is determined uniquely by the sequence $\{f_n\}_{n\in\Z}$ up
to a polynomial of degree at most $1$, and a polynomial $g$ belongs to 
$B^1_{\be,1}\big(\R^d\big)$
if and only if $g$ is constant.

Note that the functions $f_n$ have the following properties: $f_n\in L^\be(\R^d)$ and 
$\supp\F f\subset\{\xi\in\R^d:~\|\xi\|\le2^{n+1}\}$. Bounded functions whose Fourier transforms are supported in $\{\xi\in\R^d:~\|\xi\|\le\s\}$ can be characterized by the following Paley--Wiener--Schwartz type theorem  (see \cite{R}, Theorem 7.23 and exercise 15 of Chapter 7):

{\it Let $f$ be a continuous function
on $\R^d$ and let $M,\,\s>0$. The following statements are equivalent:

{\em(i)} $|f|\le M$ and $\supp\F f\subset\{\xi\in\R^d:\|\xi\|\le\s\}$;

{\em(ii)} $f$ is a restriction to $\R^d$ of an entire function on $\C^d$ such that 
$$
|f(z)|\le Me^{\s\|\im z\|}
$$
for all $z\in\C^d$.}

\

\section{\bf The main result}
\setcounter{equation}{0}

\

In this section we show that there is no Lipschitz type estimate in the trace norm for functions in the Besov space $B_{\be,1}^1(\R^3)$ and triples of not necessarily commuting self-adjoint operators of finite rank. 

For a positive number $\s$ and a positive integer $d$, we denote by $\E^\be_\s(\R^d)$ the class of functions $g$ in $L^\be(\R^d)$ such that $\supp g\subset[-\s,\s]^d$.

Recall that for a differentiable function $\psi$ on $\R$ the divided difference $\dg\f$ is defined by
$$
(\dg\f)(x,y)=\left\{\begin{array}{ll}\frac{\f(x)-\f(y)}{x-y},&x\ne y\\[.2cm]
\f'(x),&x=y.\end{array}
\right..
$$


\begin{thm}
\label{Lip3S1}
There is no constant $C>0$ such that 
\begin{align}
\label{nerC}
\|f(A_1,B_1,C_1)&-f(A_2,B_2,C_2)\|_{\bS_1}\nonumber\\[.2cm]
&\le C\|f\|_{L^\be(\R^3)}
\max\big\{\|A_1-A_2\|_{\bS_1},\|B_1-B_2\|_{\bS_1},\|C_1-C_2\|_{\bS_1}\big\}
\end{align}
for all triples of finite rank self-adjoint operators $(A_1,B_1,C_1)$ and $(A_2,B_2,C_2)$ and all functions $f$ in $\E_1^\be(\R^3)$.
\end{thm}

\Pf It is certainly enough to consider the case when $A_1=A_2=A$, $C_1=C_2=C$ and
$f(x,y,z)=\f(x,z)\psi(y)$, where $\f\in\E^\be_1(\R^2)$ and $\psi\in\E^\be_1(\R)$. 

Put $Q\df\psi(B_1)-\psi(B_2)$. We have
\begin{align}
\label{fB1-fB2}
Q&=
\iint\limits_{\R\times\R}
\big(\dg\psi\big)(y_1,y_2)\,dE_{B_1}(y_1)(B_1-B_2)\,dE_{B_2}(y_2)\nonumber\\[.2cm]
&=\sum_{\l\in\s(B_1),\,\mu\in\s(B_2),\,\l\ne\mu}\frac{\psi(\l)-\psi(\mu)}{\l-\mu}
E_{B_1}(\l)(B_1-B_2)E_{B_2}(\mu).
\end{align}
Since our operators have finite rank, the sum in 
\rf{fB1-fB2} is finite and the verification of \rf{fB1-fB2} is an elementary exercise.

We have the following formula:
\begin{align}
\label{fEQE}
f(A,B_1,C)-f(A,B_2,C)&=\iint\f(x,z)\,dE_A(x)Q\,dE_C(z)\nonumber\\[.2cm]
&=\sum_{\l\in\s(B_1),\,\mu\in\s(B_2)}\f(\l,\mu)E_{B_1}(\l)QE_{B_2}(\mu)
.
\end{align}
Again, the sum in \rf{fEQE} is finite and \rf{fEQE} can be established elementarily.

We are going to use an interpolation technique in $\E^\be_\s(\R^d)$, see \cite{ANP}.

We need the following lemma:

\begin{lem}
\label{diskSM}
Suppose that inequality {\em\rf{nerC}} holds for some positive number $C$ and for all $f$, $(A_1,B_1,C_1)$ and $(A_2,B_2,C_2)$ satisfying the hypotheses of Theorem 
{\em\ref{Lip3S1}}. Then $\E_1^\be(\R^2)\subset\fM(\R\times\R)$ and
$$
\|\f\|_{\fM(\R\times\R)}\le\const\|\f\|_{L^\be(\R^2)}.
$$
\end{lem}

\Pf Consider the function $\eta$ on $\R$ defined by
$$
\eta(x)=\frac{2(1-\cos x)}{x^2},\quad x\in\R.
$$
It is well known that $\eta\in\E_1^\be(\R)$. Clearly, $\eta(0)=1$ and $\eta(2k\pi)=0$,
$k\in\Z\setminus\{0\}$. Put $\psi(x)\df\eta(x-2\pi)$, $x\in\R$.

Let $P$ be a rank one orthogonal projection, $B_1=2\pi P$ 
and $B_2=\0$. It is easy to see that $\psi(B_1)=P$ and $\psi(B_2)=\0$.

It follows from \rf{Lip3S1} that
$$
\|f(A,2\pi P,C)-f(A,\0,C)\|_{\bS_1}
=\left\|\iint\f(x,z)\,dE_A(x)P\,dE_C(z)\right\|_{\bS_1}
\le\const\|\f\|_{L^\be(\R^2)}
$$
for arbitrary finite rank self-adjoint operators $A$ and $B$, an arbitrary rank one projection $P$ and an arbitrary function $\f$ in $\E^\be_1(\R^2)$. This implies that
$$
\left\|\iint\f(x,z)\,dE_A(x)T\,dE_C(z)\right\|_{\bS_1}
\le\const\|\f\|_{L^\be(\R^2)}\|T\|_{\bS_1}
$$
for an arbitrary function $\f$ in $\E^\be_1(\R^2)$ and an arbitrary trace class operator $T$. It is easy to see that this is equivalent to the conclusion of the lemma. $\bl$

\medskip

Let us complete the proof of Theorem \ref{Lip3S1}. Suppose that $\{c_{jk}\}$ is a finite family of complex numbers. Define the function $\f$ by
\bay
\label{oprpsi}
\f(x,y)=\sum_{j,k}\eta_j(x)\eta_k(y).
\ey
where $\eta_j(x)\df\eta(x-2\pi j)$. Then  $\f\in\E^\be_1(\R^2)$ and
\bay
\label{otscjk}
\|\f\|_{L^\be(\R^2)}\le\const\sup_{j,k}|c_{jk}|,
\ey
see \cite{ANP}, \S\;8.

It follows that we can define $\f$ by \rf{oprpsi} for an infinite family 
$\{c_{jk}\}$ of bounded complex numbers. Then the function $\f$ belongs to 
$\E^\be_1(\R^2)$ and \rf{otscjk} holds. It is also easy to see that
$\f(2\pi j,2\pi k)=c_{jk}$. 

Together with Lemma \ref{diskSM} this implies that an arbitrary bounded matrix
$\{c_{jk}\}_{j,k\ge0}$ is a Schur multiplier. However, it is well known that this is
false. In particular, the matrix $\{c_jk\}_{j,k\ge0}$ defined by
$$
c_{jk}=\left\{\begin{array}{ll}1,&j\le k\\0,&j>k,\end{array}\right.
$$
induces the operator of triangular projection, which is unbounded on $\bS_1$, see \cite{GK}. $\bl$

\medskip

Theorem \ref{Lip3S1} shows that there is no Lipschitz type estimate in the trace norm. However, in the construction given above the norm of the difference $B_1-B_2$ is separated away from zero. The following theorem shows how to make the norm of the difference as small as possible.

\begin{thm}
\label{malraz}
There exist a sequence $\{g_n\}_{n\ge0}$ of functions in $B_{\be,1}^1(\R^3)$ and sequences of self-adjoint finite rank operators $\big\{A^{(n)}\big\}_{n\ge0}$, 
$\big\{B_1^{(n)}\big\}_{n\ge0}$, $\big\{B_2^{(n)}\big\}_{n\ge0}$
and $\big\{C^{(n)}\big\}_{n\ge0}$
such that the norms $\|g_n\|_{B_{\be,1}^1}$ are bounded,
$$
\lim_{n\to\be}\big\|B_1^{(n)}-B_2^{(n)}\big\|_{\bS_1}\to0,
$$
but
$$
\big\|g_n\big(A^{(n)},B_1^{(n)},C^{(n)}\big)-
g_n\big(A^{(n)},B_2^{(n)},C^{(n)}\big)\big\|_{\bS_1}\to\be.
$$
\end{thm}

\Pf The  proof of 
Theorem \ref{Lip3S1} allows us to find sequences $f_n\in\E^\be_1(\R^3)$, finite rank
self-adjoint operators $A_n$, $B_{1,n}$, $B_{2,n}$ and $C_n$ such that 
$$
\|f_n\|_{L^\be(\R^3)}\le\const,\quad
\|B_{1,n}-B_{2,n}\|_{\bS_1}\le\const,
$$
but 
$$
\big\|f_n(A_n,B_{1,n},C_n)-
f_n(A_n,B_{2,n},C_n)\big\|_{\bS_1}\to\be.
$$

Consider a sequence $\{\e_n\}$ of positive numbers that tends to 0. Put 
$g_n(x,y,z)\df\e_nf_n\big(\frac{x}{\e_n},\frac{y}{\e_n},\frac{z}{\e_n}\big)$. Then 
$$
\|g_n\|_{B_{\be 1}^1(\R^3)}=\|f_n\|_{B_{\be 1}^1(\R^3)},
$$ 
\begin{align*}
\|g_n(\e_nA_n,\e_nB_{1,n},\e_nC_n)&-g_n(\e_nA_n,\e_nB_{2,n},\e_nC_n)\|_{\bS_1}\\[.2cm]
&=\e_n\|f_n(A_n,B_{1,n},C_n)-f(A_n,B_{2,n},C_n)\|_{\bS_1}
\end{align*}
and
$$
\|\e_nB_{1,n}-\e_nB_{2,n}\|_{\bS_1}=\e_n.
$$
It remains to choose a suitable sequence $\{\e_n\}$ and put $A^{(n)}=\e_nA_n$,
$B_1^{(n)}=\e_nB_{1,n}$, $B_2^{(n)}=\e_nB_{2,n}$ and $C^{(n)}=\e_nC_n$. $\bl$

\

\

\noindent
Department of Mathematics\\
Michigan State University\\
East Lansing Michigan 48824\\

\end{document}

%% file: classstartscr.tex
\newcommand{\lb}{\linebreak}

\newcommand{\g}{\gamma}

\newcommand{\e}{\varepsilon}

\renewcommand{\l}{\lambda}

\newcommand{\s}{\sigma}

\newcommand{\f}{\varphi}

\newcommand{\E}{{\mathscr E}}

\newcommand{\F}{{\mathscr F}}

\newcommand{\h}{{\mathscr H}}

\newcommand{\X}{{\mathscr X}}
\newcommand{\Y}{{\mathscr Y}}

\newcommand{\C}{{\Bbb C}}

\newcommand{\R}{{\Bbb R}}
\newcommand{\Z}{{\Bbb Z}}

\newcommand{\0}{{\boldsymbol{0}}}

\newcommand{\bs}{\boldsymbol}

\newcommand{\bS}{{\boldsymbol S}}

\newcommand{\rf}[1]{(\ref{#1})}

\newcommand{\df}{\stackrel{\mathrm{def}}{=}}

\newcommand{\supp}{\operatorname{supp}}

\newcommand{\const}{\operatorname{const}}

\newcommand{\eeq}{\end{equation}}
\newcommand{\beq}{\begin{equation}}
\newcommand{\bay}{\begin{eqnarray}}
\newcommand{\ba}{\begin{align*}}
\newcommand{\ea}{\end{align*}}
\newcommand{\ey}{\end{eqnarray}}
\newcommand{\bey}{\begin{eqnarray*}}
\newcommand{\eey}{\end{eqnarray*}}

\newcommand{\be}{\infty}

\newcommand{\bl}{\blacksquare}

\newcommand{\Pf}{{\bf Proof. }}
\newcommand{\im}{\operatorname{Im}}

\newtheorem{thm}{\hspace{\parindent}Theorem}[section]

\newtheorem{lem}[thm]{\hspace{\parindent}Lemma}

%% file: 3nekom.bbl
\begin{thebibliography}{99}
\label{bibl}


%

%
%
%

\bibitem[AP]{AP}{\sc A.B. Aleksandrov} and {\sc V.V. Peller}, {\it Operator Lipschitz functions}, to appear.

\bibitem[ANP]{ANP}{\sc A.B. Aleksandrov, F.L. Nazarov} and
{\sc V.V. Peller}, {\em Triple operator integrals in Schatten--von Neumann norms and functions of perturbed noncommuting operators}, C.R. Acad. Sci. Paris, S\'er. I {\bf353} (2015), 723–-728.


\bibitem[APPS]{APPS} {\sc A.B. Aleksandrov, V.V. Peller, D. Potapov}, and
{\sc F. Sukochev}, {\em Functions of normal operators under perturbations},
Advances in Math. {\bf226} (2011), 5216-–5251.


\bibitem[BS1]{BS1} {\sc M.S. Birman} and {\sc M.Z. Solomyak},
{\em Double Stieltjes operator integrals},
Problems of Math. Phys., Leningrad. Univ. {\bf1} (1966), 33--67 (Russian).
English transl.: Topics Math. Physics {\bf1} (1967), 25--54, Consultants Bureau Plenum
Publishing Corporation, New York.

\bibitem[BS2]{BS2} {\sc M.S. Birman} and {\sc M.Z. Solomyak},
 {\em Double Stieltjes operator integrals. II},
 Problems of Math. Phys., Leningrad. Univ. {\bf2} (1967), 26--60 (Russian).
English transl.: Topics Math. Physics {\bf2} (1968), 19--46, Consultants Bureau Plenum
Publishing Corporation, New York.

%

\bibitem[BS3]{BS3} {\sc M.S. Birman} and {\sc M.Z. Solomyak},
{\em Double Stieltjes operator integrals. III},
Problems of Math. Phys., Leningrad. Univ. {\bf6} (1973), 27--53 (Russian).




\bibitem[DK]{DK} {\sc Yu.L. Daletskii} and {\sc S.G. Krein}, {\em Integration and differentiation of
functions of Hermitian operators and application to the theory of perturbations} (Russian), Trudy Sem.
Functsion. Anal., Voronezh. Gos. Univ. {\bf1} (1956), 81--105.


\bibitem[F1]{F1}  {\sc Yu.B. Farforovskaya}, {\em  The connection of the Kantorovich-Rubinshtein metric for spectral resolutions of selfadjoint operators with functions of operators},
Vestnik Leningrad. Univ.  {\bf19}  (1968), 94--97. (Russian).

\bibitem[F2]{F2}  {\sc Yu.B. Farforovskaya}, {\em An example of a Lipschitzian function of selfadjoint
operators that yields a nonnuclear increase under a nuclear perturbation}.  Zap. Nauchn. Sem.
Leningrad. Otdel. Mat. Inst. Steklov. (LOMI)  {\bf30}  (1972), 146--153 (Russian).



%

%

\bibitem[GK]{GK} {\sc I.C. Gohberg, M.G. Krein}, {\em Theory of Volterra operators in Hilbert space and its
applications}, Nauka, Moscow, 1965; English translation: American Mathematical Society, Providence, R.I. 1970.

\bibitem[JTT]{JTT} {\sc K. Juschenko, I.G. Todorov} and {\sc L. Turowska}, {\em Multidimensional operator multipliers}, Trans. Amer. Math. Soc. {\bf361}
(2009), 4683-–4720.

%

\bibitem[Ka]{Ka}  {\sc T. Kato}, {\em Continuity of the map $S\mapsto \mid S\mid $ for linear operators},
Proc. Japan Acad.  {\bf49}  (1973), 157--160.


\bibitem[KPSS]{KPSS} {\sc E. Kissin, D. Potapov, V. S. Shulman} and {\sc F. Sukochev}, {\it Operator smoothness in Schatten norms for functions of several variables: Lipschitz conditions, differentiability and unbounded derivations},  Proc. Lond. Math. Soc. (3) {\bf105} (2012), 661--702.





\bibitem[Mc]{Mc}
{\sc A. McIntosh}, {\em Counterexample to a question on commutators},
Proc. Amer. Math. Soc. {\bf29} (1971) 337--340.

\bibitem[NP]{NP} {\sc F.L. Nazarov} and {V.V. Peller} {\em Functions of $n$-tuples of commuting self-adjoint operators}, J. Funct. Anal. {\bf266} (2014), 5398–-5428.


\bibitem[Pee]{Pee} {\sc J. Peetre},
{\em New thoughts on Besov spaces}, Duke Univ. Press., Durham, NC, 1976.


\bibitem[Pe1]{Pe1} {\sc V.V.Peller}, {\em Hankel operators of class ${\bf S}_{p}$
and their applications (rational approximation, Gaussian processes,
the problem of majorizing operators)}, Mat. Sbornik,
{\bf 113} (1980), 538-581.

English Transl. in Math. USSR Sbornik, {\bf 41}
(1982), 443-479.


\bibitem[Pe2]{Pe2} {\sc V.V. Peller},
{\em Hankel operators in the theory of perturbations of unitary and self-adjoint operators},
Funktsional. Anal. i Prilozhen. {\bf19:2}  (1985),
37--51 (Russian).
English transl.: Funct. Anal. Appl. {\bf19} (1985) , 111--123.



\bibitem[Pe3]{Pe3} {\sc V.V. Peller},
{\em Hankel operators in the perturbation theory of of unbounded self-adjoint operators}.
Analysis and partial differential equations,  529--544,
Lecture Notes in Pure and Appl. Math., {\bf122}, Dekker, New York, 1990.

\bibitem[Pe4]{Pe4} {\sc V.V. Peller}, {\em Hankel operators and their applications,}
Springer-Verlag, New York, 2003.

\bibitem[Pe5]{Pe5} {\sc V.V. Peller}, {\em Multiple operator integrals and higher operator
derivatives}, J. Funct. Anal.  {\bf233}  (2006),  515--544.

\bibitem[Pe6]{Pe6} {\sc V.V. Peller}, {\em Multiple operator integrals in perturbation theory}, Bull. Math. Sci. {\bf6} (2016), 15--88.


\bibitem[Pi]{Pi} {\sc G. Pisier}, {\em Similarity problems and completely bounded maps},
Second, expanded edition. Includes the solution to ``The Halmos problem''. Lecture Notes in Mathematics,
1618. Springer-Verlag, Berlin, 2001.

\bibitem[PS]{PS} {\sc D. Potapov} and {\sc F. Sukochev}, {\em Operator-Lipschitz functions in Schatten--von Neumann classes}, Acta Math. {\bf207} (2011), 375--389.

\bibitem[R]{R}{\sc W. Rudin}, {\em Functional analysis}, M$^{\rm c}$Graw Hill, 1991.


%

%









%
%

\end{thebibliography}
